\documentclass[12pt,a4paper]{article}
\usepackage[margin=27mm]{geometry}
\usepackage{microtype,hyperref,xcolor,titlesec,fancyhdr}
\usepackage{xurl}
\usepackage{listings}
\usepackage{iftex}
\newtheorem{theorem}{Theorem}
\ifPDFTeX
  \DeclareUnicodeCharacter{03C4}{\ensuremath{\tau}}
\else
  \usepackage{fontspec}
\fi
\hypersetup{colorlinks=true,linkcolor=black,urlcolor=blue!55!black,citecolor=blue!55!black,pdftitle={Moonshine Conjecture Formalized},pdfauthor={Research Team}}
\makeatletter
\renewcommand*\l@section{\@dottedtocline{1}{1.5em}{2.3em}}
\makeatother
\titleformat{\section}{\fontsize{14}{17}\selectfont\bfseries}{\thesection}{0.7em}{}
\titleformat{\subsection}{\fontsize{13}{16}\selectfont\bfseries}{\thesubsection}{0.6em}{}
\titlespacing*{\section}{0pt}{2.5ex plus .5ex minus .2ex}{1.2ex}
\titlespacing*{\subsection}{0pt}{2.0ex plus .4ex minus .2ex}{0.9ex}
\newcommand{\code}[1]{{\fontsize{10}{12}\selectfont\ttfamily #1}}
\definecolor{leanKeyword}{HTML}{1F5F8B}
\definecolor{leanSymbol}{HTML}{7A3E9D}
\lstdefinelanguage{Lean}{morekeywords={theorem,fun},sensitive=true,morecomment=[l]{--},morestring=[b]"}
\lstdefinestyle{leanCode}{
  language=Lean,
  basicstyle=\ttfamily\fontsize{10}{12}\selectfont,
  keywordstyle=\color{leanKeyword}\bfseries,
  commentstyle=\color{gray!70!black}\itshape,
  columns=fullflexible,
  keepspaces=true,
  showstringspaces=false,
  frame=single,
  framerule=0.3pt,
  rulecolor=\color{black!20},
  backgroundcolor=\color{black!3},
  xleftmargin=0.4em,
  xrightmargin=0.4em,
  aboveskip=0.8em,
  belowskip=0.8em,
  literate=
    {τ}{{\textcolor{leanSymbol}{\ensuremath{\tau}}}}1
    {∧}{{\textcolor{leanSymbol}{\ensuremath{\wedge}}}}1
    {∃}{{\textcolor{leanSymbol}{\ensuremath{\exists}}}}1
    {≃}{{\textcolor{leanSymbol}{\ensuremath{\simeq}}}}1
}
\begin{document}
\begin{center}
{\LARGE\bfseries Moonshine Conjecture Formalized\par}
\vspace{0.45em}
{\large Toward an Agentic Era of Lean Formalization\par}
\vspace{1em}
{\large Jin-Cheng Guu\footnote{University of Alberta. \href{mailto:guu@ualberta.ca}{guu@ualberta.ca}}\par}
\end{center}
\begin{abstract}
The celebrated and intimidating Moonshine Conjecture has now been formalized
in Lean \cite{MoonshineMinimal2026,VOARepo1062026}. This achievement might suggest
that, once a sufficiently large library compiles, the principal mathematical
work is complete. This paper tells a less optimistic story.

During our formalization project, we repeatedly encountered two problems that
a formally checked proof alone does not settle: (1) whether the formal
statements actually express the mathematics we humans expect, and (2) how
that meaning survives as the repositories and their dependencies evolve.

We argue that there is no magic pill that turns a large agent-produced library
into mathematics a community can trust. Even after the proofs are checked by the
Lean kernel, someone must judge whether the definitions and statements express
the intended mathematics and whether future updates preserve that meaning. This
work, even when done, can only be trusted. We must therefore return to the usual
trust model: relying on the mathematical judgment of maintainers and expert
communities.

But semantic review is expensive, making it tempting to keep producing code and
take successful kernel checks as sufficient grounds for trust. This shortcut
accumulates semantic debt: definitions and results with unchecked mathematical
interpretations are reused as foundations for further work. A library may keep
compiling and produce the desired formal theorem while users lose track of which
mathematical claims they can responsibly rely on.

Reliable reuse is also a prerequisite for Terence Tao’s vision of formal
mathematics as something people can navigate, explain, learn from, and build
with \cite{Tao2026}.
This paper takes our work connecting the VOA
and CFSG repositories for the Moonshine Conjecture as a case study of the
challenges of reliable reuse as formal developments scale.
\end{abstract}

\vspace{1.2em}
\begingroup
\small
\setlength{\parskip}{0pt}
\renewcommand{\baselinestretch}{0.88}\selectfont
\tableofcontents
\endgroup
\clearpage
\section{Background and Aim}

Mathematical proofs are important not only because they guarantee correctness.
They also turn observed patterns into general principles, expose the
conditions on which conclusions depend, and create abstraction layers that let
us reason from completely new points of view.

\emph{Reading and checking} proofs, however, are extremely expensive. Most
proofs omit steps that trained readers can supply; for long arguments spanning
hundreds of papers and several fields, the details may only be accessible to a
few experts in the world. Given the nature of proofs, and because results are
reused, a single uncorrected error in a foundational argument can propagate
through later work and undermine an entire body of results. For this reason,
historically and today, human mathematicians do their best to check all
details manually. Fortunately, any proof can be broken down into many minute
steps, each of which is routine enough for a computer to check.

This pressure, together with the possibility of mechanical checking, led to the
development of formal verifiers like Lean. Once a proof is written in Lean, the
kernel mechanically checks every formal inference step against the given
definitions and assumptions. Programs like this have been developed for decades,
but writing and understanding fully detailed formal proofs have often been even
harder than writing ordinary mathematical proofs.

Fortunately, in 2026, agentic LLMs have become powerful enough to effectively
assist with these tasks. Initial large-scale successes include the formalization
of the Classification of the Finite Simple Groups, the Stacks Project, and
Fermat's Last Theorem \cite{AxiomCFSG2026,AxiomStacks2026,AnthropicFLT2026}.
These projects demonstrate both the scale agents can now reach and the review
burden that follows: their formal artifacts are too large for any one person to
read, understand, and check end to end without substantial collective effort.

This opportunity is exciting, but proof production alone does not give us formal
mathematics that a community can trust and reuse. This paper examines the limits
of current success claims and the responsibilities that follow for the
community.

\section{Semantic Gap: What Did We Formalize?}\label{sec:semantic-gap}

\noindent\textit{Are we proving the intended theorem?}\par

Lean's kernel can reliably check an enormous number of formal inference steps,
but it does so only relative to the definitions and assumptions it is given.
We still have to verify that the symbols and statements carry the meanings we intend.

For example, suppose a hypothetical library defines \code{IsPrime(p)} to mean
that \(p\) is prime in the usual sense and \(p>3\). It can then correctly prove
that every \code{IsPrime} is congruent to \(1\) or \(5\) modulo \(6\). A reader
who interprets the name as ordinary primality overlooks the counterexamples
\(2\) and \(3\). And this misunderstanding can propagate silently in the ocean
of lean codes. While this is just a toy example, in a large network of formal
libraries, similar restrictions can be harder to notice and can affect much
less familiar mathematical objects.

A naive response is to inspect every symbol in each statement one cares about.
But as one cares about more statements and more of their proofs, the burden
grows rapidly. Eventually one relies on libraries too large to read end to end;
even a complete semantic review of one version would not remove the need to
review later updates.

Recent large public releases of AI-generated mathematical
results make this scaling problem immediate: even when many arguments are
accompanied by Lean formalizations, readers cannot reasonably be expected to
independently reconstruct every definition, dependency, and intended scope
\cite{OpenAIMathRelease2026}.

\section{Versioning: Maintaining Trust}\label{sec:versioning}

\noindent\textit{Can meaning survive independent updates?}\par

Reviewing one version of a library does not settle whether later versions still
express the intended mathematics. An update may change definitions, theorem
statements, or dependencies, including those used by downstream projects. So the
users have to rely on trusted maintainers to review these changes and explain
their consequences; they cannot continually re-audit the entire ecosystem
themselves. An error correction requires more than a successful rebuild:
downstream users must know which results it affects.

This problem is familiar in software engineering but easy for mathematicians
to underestimate. A library may appear stable as a mathematical object, yet its
code and interfaces are embedded in an evolving dependency graph. Our
case in Section~\ref{sec:moonshine-case} provides a concrete example
involving incompatible Lean and Mathlib versions.

When a widely used and trusted library releases a new version, projects that
depend on it face a choice: upgrade or remain pinned to an older release. Neither
choice is costless.

\begin{itemize}
\item Upgrading can require compatibility repairs and new bridge proofs.
  Moreover, a large upstream change can impose a costly review burden on
  downstream projects. Their maintainers must determine not only whether their
  code still compiles, but whether the changed definitions and interfaces
  preserve the meanings on which their own claims depend.
\item Remaining pinned carries a different cost: If the newer release corrects
  an error, downstream artifacts that stay on the old version do not
  automatically inherit the trust claim. Even when this poses no immediate
  danger, it creates a growing burden of explanation: the project's trust claim
  may gradually erode, and catching up becomes increasingly difficult once
  versions have diverged too far.
\end{itemize}

Versioning concerns not only which definitions and dependencies a result uses,
\textbf{but also whether a community can still trust their mathematical meaning
  as those versions change}. Successful checks by Lean's kernel alone cannot
settle this issue, because they do not establish that the formal statements
express the intended mathematics. Therefore maintaining that trust across
updates takes continuing work. Maintainers need enough expertise and reputation
to review changes, explain their mathematical consequences, and take
responsibility for the versions they support.

\section{Moonshine: A Case Study}\label{sec:moonshine-case}

The Moonshine Conjecture begins with an unexpected
connection between the \(j\)-invariant, which classifies complex elliptic
curves (equivalently, complex tori), and the Monster, the largest sporadic
finite simple group. Moonshine has several formulations and extensions. In
this paper we formalize one of its most famous and illuminating
statements:
\begin{theorem}
The vertex operator algebra \(V^\natural\) has graded character
\(j(\tau)-744\), and its full automorphism group is isomorphic to the Monster.
\end{theorem}
This explains in particular the famous equation:

\[ 196884 = 196883 + 1. \]

Frenkel, Lepowsky, and Meurman constructed \(V^\natural\), giving this
correspondence an algebraic foundation. Building on their work, Borcherds
proved the Conway--Norton moonshine conjecture in 1992 using methods inspired
by string theory, an achievement recognized by a Fields Medal in 1998
\cite{GannonMoonshine2004,QuantaMoonshine2015}.

Formalizing this result required connecting two large
formal worlds: the VOA repository \cite{VOARepo2026} and AxiomMath's CFSG
repository \cite{CFSGRepo2026}.\footnote{Most of the heavy lifting was done in these two repositories.} Our
formalization \cite{MoonshineMinimal2026} proves the following
theorem\footnote{During the writing of this paper, the VOA repository released
v1.06, independently formalizing the result displayed here and a qualitative
Conway--Norton statement \cite{VOARepo1062026}. The latter proves that each
McKay--Thompson series satisfies eight analytic conditions defining a normalized
Hauptmodul. Chapter~16 of their accompanying textbook proves the equivalence
with the usual geometric definition on paper; this equivalence is not
formalized in Lean.
The formal statement asserts that
suitable groups exist; it does not identify them with the groups specified by
Conway and Norton for each Monster conjugacy class, or establish the prescribed
relation between level and element order.}:

\begin{lstlisting}[style=leanCode]
theorem moonshineConjecture :
    ∃ e : VOA.Aut Vnat ≃* MonsterGroup,
      char Vnat = fun τ => VOA.jInvariant τ - 744
\end{lstlisting}

\subsection*{Semantic Gap: What Did We Formalize?}

At first sight, this looks like a faithful translation of the mathematical
statement. But what, exactly, do its symbols mean? What is \code{Vnat}?
What is \code{jInvariant}, and is it defined as the
intended \(j\)-function in the mathematical literature? What is \code{VOA.Aut};
does it use the usual notion of a VOA automorphism, and does each automorphism
respect the conformal vector? What is \code{MonsterGroup}, and does it denote
the intended Monster group?\footnote{The Lean formalization of Fermat's Last
Theorem \cite{AnthropicFLT2026} is a useful contrast. Its statement only
involves positive integers, addition, equality, and exponentiation, all of which
can be inspected quickly.}

Indeed, these symbols come from the VOA repository and the CFSG repository. One
must inspect those repositories to check their meanings. But their dependency
graphs are too large to inspect quickly. Once one also wants to understand the
proof and its intermediate lemmas, the same questions arise for the terms
appearing there. The burden
of semantic checking can quickly grow.

\subsection*{Building a Bridge Between Repositories}

To finalize our formalization of the
Moonshine Conjecture, we had to connect the relevant VOA and CFSG results.
This required a careful bridge between the two repositories. But the CFSG
repository was developed against newer Lean and Mathlib versions than those
used by the VOA repository, so we first migrated the VOA development to
compatible dependency versions. Only after that could we construct the
required bridge.

The two repositories represented one of the key objects (the Griess algebra)
by different Lean objects, using different coordinates, bases,
quotient types, and representatives. We had to construct and verify a
coordinate equivalence, prove that it preserved the product and bilinear form,
and transport the relevant generated-subgroup statement across it.

\subsection*{Versioning: Maintaining Trust}

We could resolve these issues in our project. Our concern is how they will
scale as more independently maintained repositories are connected:

\begin{enumerate}
  \item Different repositories may use different versions of Lean, Mathlib,
        and shared libraries. Connecting them can require migrations and
        renewed semantic review. (A successful kernel check is not enough:
        maintainers must also check that the migrated definitions and
        statements preserve the intended mathematical meaning.) As more
        repositories are linked, a change in one dependency can create work
        across several projects.
  \item Our proof uses a pinned fork of the VOA repository. Lean's kernel
        checks the proof against that fork, but this does not settle who will
        maintain it, follow upstream changes, or review the versions that
        later users rely on.
\end{enumerate}

A platform that brings together mathematical repositories must also check the
meaning of the connections between them. Different formal objects may
represent the same mathematical structure, but proving an isomorphism between
them does not by itself establish that the bridge expresses the intended
mathematical identification. Its meaning depends on the objects being
connected and on what downstream results require it to preserve. That
interpretation must also survive changes on either side.

\section{Conclusions: What Remains After Proof Abundance}\label{sec:repo-obligations}

The Moonshine Conjecture has been formalized. So what? Why do we still feel
that there are gaps everywhere, and that there are endlessly many gaps to
fill?

These are not gaps in the Lean proof itself. They are all the things that
remain hard to know after the proof compiles: whether the definitions mean what
we think they mean, whether an intermediate lemma is using the right object,
whether a bridge really connects the two mathematical structures that it claims
to connect, and whether any of this remains true after the surrounding
repositories change. Once someone points to one of these questions, it may be
possible to check. The problem is that there are too many of them, spread
across too many definitions, lemmas, repositories, and versions.

This leaves us with the following lessons for building durable, reusable, and
trusted public mathematical repositories in the agentic era.

\subsection{The Bottleneck Is the Maintainer, Not the Prover}\label{sec:semantic-maintainers}

We now have more proofs, and proofs are becoming cheap. What is the next
bottleneck? As we have argued, whether semantic work has been done well cannot
itself be cheaply checked or reproduced by every later user. We are therefore
pushed back toward a familiar trust model: each repository will need a trusted
group of maintainers, and the community will have to trust that they have done
the hard work of guarding its semantics. As we explained in Section~\ref{sec:versioning},
this is not a once-and-for-all task. It requires ongoing work. This also
challenges the view that the answer is simply to collect large numbers of
contributions. A quick inspection reveals quality problems in some (if not many) existing
repositories. If repositories with high quality
standards already require sustained expert maintenance, it is unclear what a
community can responsibly do with artifacts produced without it.

There may be platforms
that collect and curate trusted combinations, but they will not make the
underlying work disappear: updates, forks, and new bridges can introduce new
semantic review obligations.

\subsection{A Version Bump Is Not a Semantic Repair}\label{sec:semantic-versioning}

\begingroup
\widowpenalty=10000
Versioning invites a false but tempting picture: find a problem, bump a
dependency, repair the errors reported by Lean, release a new version, and move on. For
a small and isolated bug, this may be enough. Whether it is enough depends on
the scope of the change. Each
correction can change definitions, interfaces, or bridges, and can therefore
change what downstream results are really using.
No package manager or successful build can
decide whether these changes still preserve the mathematical meaning that users
relied on. Maintainers earn trust not by keeping proofs accepted by Lean's kernel, but by
doing the harder work of understanding these changes, explaining their
consequences, and deciding what their community should still stand behind.
\par\endgroup

\subsection{Reliable Reuse Comes Before Public Mathematics}\label{sec:tao-horizon}

Terence Tao has articulated a larger and attractive vision: formal mathematics
should not merely certify proofs, but become something that people can
navigate, explain, learn from, and build with \cite{Tao2026}. We share this
ambition. Before a formal artifact can serve as trusted public mathematics, it must first
be something that a community can reliably reuse: its meanings must be
understood well enough, and its versions maintained carefully enough, that
later readers have something stable to explain, teach, or extend.

\appendix
\section*{Acknowledgements}
\addcontentsline{toc}{section}{Acknowledgements}

We thank Fedor Popov, Adam Topaz, and Terry Gannon for discussions.

{\fontsize{8}{10}\selectfont
\renewcommand{\code}[1]{\texttt{#1}}

}


\begin{thebibliography}{VOA1.06}
\setlength{\itemsep}{0pt}
\setlength{\parskip}{0pt}
\bibitem[OAI26]{OpenAIMathRelease2026} OpenAI, \emph{Sharing AI progress in mathematics}, 6 October 2026, \url{https://openai.com/index/sharing-ai-progress-in-mathematics/}.
\bibitem[Tao26]{Tao2026} Terence Tao, \emph{Mathematics in the Age of AI}, arXiv:2608.16753, 2026, \url{https://arxiv.org/abs/2608.16753}.
\bibitem[Ax26a]{AxiomCFSG2026} Axiom, \emph{Formalizing the Enormous Theorem}, 26 October 2026, \url{https://axiommath.ai/research/the-enormous-theorem-formalized/}.
\bibitem[CFSG]{CFSGRepo2026} AxiomMath, \emph{CFSG}, source commit \code{69cb972f}, \href{https://github.com/AxiomMath/CFSG/tree/69cb972f587a3ecc035e092253167ff2e303c962}{GitHub snapshot}.
\bibitem[Ax26b]{AxiomStacks2026} Axiom, \emph{The Stacks Project, Formalized in Lean}, 4 November 2026, \url{https://axiommath.ai/research/the-stacks-project-formalized/}.
\bibitem[Anth26]{AnthropicFLT2026} Anthropic, \emph{Formalizing Fermat's Last Theorem}, 2026, \url{https://www.anthropic.com/research/formalizing-fermats-last-theorem}.
\bibitem[Guu26]{MoonshineMinimal2026} J. C. Guu et al., \emph{Moonshine Conjecture formalization}, release v0.0.2, source commit \code{5b75562}, \href{https://github.com/jcguu95/moonshine-conjecture/tree/5b755622583b1136666049800c9bea4ead7492ba}{GitHub snapshot}.
\bibitem[VOA1.05]{VOARepo2026} VOA repository, \emph{Lean\_VOA}, v1.05, source commit \code{013a799}, \href{https://github.com/Fedor-Popov/Lean_VOA/tree/013a799ca834a8a6b3afa8129e3e522025b3f45c}{GitHub snapshot}.
\bibitem[VOA1.06]{VOARepo1062026} VOA repository, \emph{Lean\_VOA}, v1.06, source commit \code{799c7a9}, \href{https://github.com/Fedor-Popov/Lean_VOA/tree/799c7a9f3eee9226beb9a567c6bfe0fd275cb8ed}{GitHub snapshot}.
\bibitem[Mathlib]{Mathlib2026} The mathlib community, \emph{mathlib}, source commit \code{b63f6e8}, \href{https://github.com/leanprover-community/mathlib4/tree/b63f6e8a68d220e3b3bc4f3792bb53650d375f24}{GitHub snapshot}.
\bibitem[Kla15]{QuantaMoonshine2015} Erica Klarreich, \emph{Mathematicians Chase Moonshine's Shadow}, Quanta Magazine, 12 March 2015, \url{https://www.quantamagazine.org/mathematicians-chase-moonshine-string-theory-connections-20150312/}.
\bibitem[Gan04]{GannonMoonshine2004} Terry Gannon, \emph{Monstrous Moonshine: The first twenty-five years}, arXiv:math/0402345, 2004, \url{https://arxiv.org/abs/math/0402345}.
\end{thebibliography}
\end{document}